\documentclass[12pt,reqno]{article}
\usepackage[usenames]{color}
\usepackage[colorlinks=true,
linkcolor=webgreen,
filecolor=webbrown,
citecolor=webgreen]{hyperref}
\definecolor{webgreen}{rgb}{0,.5,0}
\definecolor{webbrown}{rgb}{.6,0,0}

\usepackage{epsfig}                                                                               
\usepackage{amsmath,amssymb,eucal}
\usepackage{amsfonts}
\usepackage{verbatim} 
\usepackage{multirow}

\newcommand{\rtj}[1]{\makebox[0.8cm][r]{#1}}
\newcommand{\rtjl}[1]{\makebox[1.0cm][r]{#1}}
\newcommand{\myemptyset}{\varnothing}

\newcommand{\forward}{F}
\newcommand{\backward}{B}
\newcommand{\winning}{W}
\newcommand{\boardstart}{b_s}
\newcommand{\boardend}{b_f}
\newcommand{\xstart}{x_s}
\newcommand{\ystart}{y_s}
\newcommand{\xend}{x_f}
\newcommand{\yend}{y_f}
\newcommand{\descendants}{\mathcal{D}}
\newcommand{\super}{\mathfrak}
\newcommand{\example}{F}

\newenvironment{packed_enumerate}{
\setlength{\parsep}{0pt}
\setlength{\parskip}{0pt}
\begin{enumerate}
  \setlength{\itemsep}{1pt}
  \setlength{\parsep}{0pt}
  \setlength{\parskip}{0pt}
}{\end{enumerate}}

\begin{document}

\begin{center}
\vskip 1cm{\LARGE\bf 
Notes on solving and playing\\peg solitaire on a computer}
\vskip 1cm
{\large George I. Bell}\\
\href{mailto:gibell@comcast.net}{\tt gibell@comcast.net} \\
\end{center}

\vskip .2 in
\begin{abstract}
We consider the one-person game of peg solitaire played on a computer.
Two popular board shapes are the 33-hole cross-shaped board,
and the 15-hole triangle board---we use them as examples throughout.
The basic game begins from a full board with one peg missing
and the goal is to finish at a board position with one peg.
First, we discuss ways to solve the basic game on a computer.
Then we consider the problem of quickly distinguishing board positions
where the goal can still be reached (``winning" board positions) from those where it cannot.
This enables a computer to alert the player if a jump under consideration leads to a dead end.
On the 15-hole triangle board, it is possible to identify all winning board positions
(from any single vacancy start) by storing a key set of 437 board positions.
For the ``central game" on the 33-hole cross-shaped board,
we can identify all winning board positions by storing 839,536 board positions.
By viewing a successful game as a traversal of a directed graph of winning board positions,
we apply a simple algorithm to count the number of ways to traverse this graph,
and calculate that the total number of solutions to the central game is
40,861,647,040,079,968.
Our analysis can also determine how quickly we can reach a ``dead board position",
where a one peg finish is no longer possible.
\end{abstract}

\section{Introduction}
\label{sec:intro}

\begin{figure}[htb]
\centering
\epsfig{file=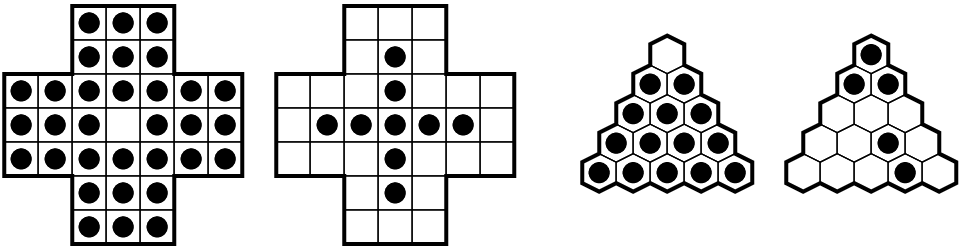}
\caption{Popular peg solitaire boards: (a,b) the 33-hole cross-shaped board,
(c,d) the 15-hole triangle board.}
\label{fig1}
\end{figure}

Peg solitaire is one of the oldest puzzles,
with a 300 year history.
The puzzle consists of a game board together with a number of pegs,
or more commonly marbles.
The board contains a grid of holes in which these pegs or marbles are placed.
Figure~\ref{fig1} shows the most common shapes for a peg solitaire board,
the 33-hole cross-shaped board, and the 15-hole triangle board.

The game is played by jumping one peg over another
into an empty hole,
removing the peg that was jumped over.
On the cross-shaped board of Figure~\ref{fig1}a,
these jumps must be made along columns or rows,
whereas on the triangle board of Figure~\ref{fig1}c,
jumps are allowed along any of the six directions parallel
to the sides of the board.
The goal is to finish with one peg,
a more advanced variation is to finish with one peg
at a specified hole.
The basic game begins from a full board with one
peg removed, as in Figure~\ref{fig1}a or c.
Starting from Figure~\ref{fig1}a and ending with one peg in
the center of the board is known as the ``central game" \cite{Beasley}.
In Figure~\ref{fig1}b or d we show starting configurations to be
referred to later.

Only relatively recently in the history of the game have computers
been used as an interface to play the game,
as well as solve the game.
There are now dozens of versions of the game available for playing
on your computer or even your cell phone.

A computer version of the game is in many ways less satisfying
than a physical game.
However, there are some definite advantages to playing
the game on a computer.
The board can be reset instantly,
and you won't be chasing marbles that fall off the board!
You can take back a jump,
all the way back to the beginning if desired.
This tends to make the game easier as you can more easily backtrack
from dead ends.
The sequence of jumps leading to a solution
can be recorded and played back.
The computer can also be programmed
to tell the user if the jump they are considering leads
to a dead end or not.
Adding this ability to a computer version of the game is tricky,
and most versions do not have this ability.
It is the goal of this paper to describe efficient techniques
to enable a computer to point out all good and bad jumps from the
current board position.

\section{Board types and symmetry}
\label{sec:boards}

\begin{figure}[htb]
\centering
\epsfig{file=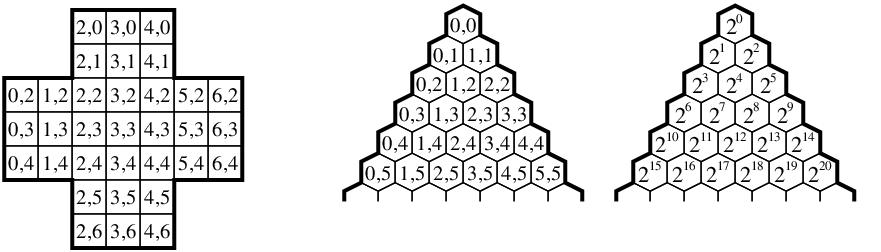}
\caption{Hole coordinates for (a) the 33-hole cross-shaped board, (b) the triangle boards.
(c) The weighting of each hole to convert a board position to a binary code.}
\label{fig2}
\end{figure}

We label the holes in the 33-hole board using Cartesian coordinates
(Figure~\ref{fig2}a), but with $y$ \textit{increasing} downward.
For the triangle boards, we use ``skew-coordinates" as shown
in Figure~\ref{fig2}b.
By adding 1 to each coordinate, and converting the first
to a letter, we obtain the standard board labelings used by
Beasley \cite{Beasley} and Bell \cite{BellSol}
(for example the central hole $(3,3)$ in Figure~\ref{fig2}a becomes ``d4",
see Figure~\ref{fig_notation}).

The 33-hole board has square symmetry.
There are eight symmetry transformations of the board,
given by the identity, rotations of $90^\circ$, $180^\circ$ and $270^\circ$,
and a reflection of the board followed by these 4 rotations \cite{BellFr}
(the dihedral group $D_4$).
The triangle boards have 6-fold symmetry, with 3 possible
rotations of $0^\circ$, $120^\circ$ or $240^\circ$, plus a reflection
followed by a rotation (the dihedral group $D_3$).

To store a particular board position on a computer,
we convert it to an integer by taking one bit per hole.
The most obvious way to do this is to take the board, row by row,
top to bottom, as in Figure~\ref{fig2}c.
We will use $N$ to denote the total number of holes on the board (the board size),
so each board position is represented by an $N$ bit integer.
If $b$ is a board position we'll denote this integer representation by
$\mbox{code}(b)$.
Many computer languages use a 4-byte integer,
the 33-hole board needs one more bit!
Beasley \cite[p. 249]{Beasley} gives a technique for storing a board position
on the 33-hole board using 4 fewer bits.
This technique is used in the
\href{http://arxiv.org/src/0903.3696v4/anc}{ancillary program}
``pegs.cpp".
For boards with more than 32 holes, we usually split
$\mbox{code}(b)$
into several 4-byte integers.

The \textbf{complement} of a board position $b$ is obtained by
replacing every peg by a hole (i.e. removing it), and replacing
every hole by a peg.
The complement of $b$ will be denoted as $\overline{b}$.
We note that $\mbox{code}(\overline{b})=\mbox{code}(f)-\mbox{code}(b)$,
where $f$ is the board position where every hole contains a peg,
$\mbox{code}(f)=2^N-1$.
The starting position for the ``central game" in Figure~\ref{fig1}a
therefore has code $2^{33}-2^{16}-1$.

Two board positions are symmetry equivalent if one can be converted to
the other by a symmetry transformation.
This equivalence relation introduces a set of \textit{equivalence classes}
of board positions, which we call \textbf{symmetry classes}.
The symmetry class of a board position does not change after it
is rotated or reflected.
One way to choose a representative from each symmetry class is
to take the one with the smallest code.
We use the notation $\mbox{mincode}(b)$ to denote this operation.
For example, the board position $b$ in Figure~\ref{fig1}d has code
$2^0+2^1+2^2+2^{8}+2^{13}=8455$,
and the other 5 codes obtained by symmetry transformation are:
$2183$, $3156$, $3904$, $25106$ and $25280$,
so $\mbox{mincode}(b)=2183$.
We also have $\mbox{mincode}(\overline{b})=(2^{15}-1)-\mbox{maxcode}(b)=7487$.

\section{Single vacancy to single survivor problems}
\label{sec:SVSS}

A peg solitaire problem which begins with one peg missing,
with the goal to finish with one peg, will be called
a \textbf{single vacancy to single survivor} problem,
abbreviated SVSS.
When the starting vacancy $(\xstart,\ystart)$ and finishing hole
$(\xend,\yend)$ are the same, the SVSS problem
is called a \textbf{complement problem},
because the starting and ending board positions
are complements of one another.

A simple parity argument gives a necessary condition for
solvability of a SVSS problem \cite[Chapter 4]{Beasley}.
On a square lattice (like the standard 33-hole board),
the requirement is that $\xstart$ and $\xend$ must differ by a multiple of 3,
or that $\xstart\equiv \xend \pmod 3$, and $\ystart\equiv \yend \pmod 3$.
Starting and ending board positions satisfying the above conditions
are said to be in the same \textbf{position class}
(a second equivalence class).
On a triangular lattice, the requirement is weaker:
$\xstart+\ystart \equiv \xend+\yend \pmod 3$\footnote{Both
these conditions assume that the full and empty boards are in the same
position class, a board satisfying this is called a \textbf{null-class} board
(see Beasley \cite{Beasley}).
All the boards we will consider are
null-class, except for the triangle board of side $4$.}.
The main result is
\textit{a peg solitaire jump does not change the position class}, 
thus an entire peg solitaire game is played in the same position class.
We will not go into the theory of position classes,
the reader should see Beasley \cite[Chapter 4]{Beasley} or
Bell \cite{BellSol} for triangular peg solitaire.

It is interesting to see what happens to the position
class after the board is rotated or reflected.
For the central game on the 33-hole board (Figure~\ref{fig1}a),
the position class is not changed by rotations or reflections
of the board.
If we begin with $(3,3)$ vacant we can finish at $(3,3)$,
or the rotationally equivalent holes $(3,0)$, $(0,3)$, $(6,3)$ and $(3,6)$.
Any board position which begins from any of these five holes
is in the position class of one peg in the center.
Moreover, if we reflect and/or rotate the board at any time
during the game, it remains in the same position class.
We call any SVSS problem beginning and ending at any of these holes ``Type A";
their solutions are all interconnected (share the same board positions).

\begin{figure}[htb]
\centering
\epsfig{file=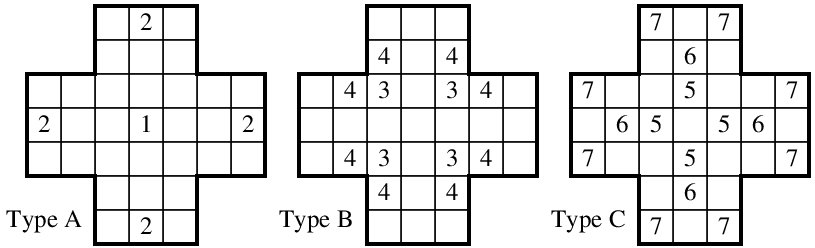}
\caption{The three types of SVSS problem on the standard 33-hole board.}
\label{fig3}
\end{figure}

More commonly, the position class does change after the board is
rotated or reflected.
``Type B" problems are shown in Figure~\ref{fig3}b.
In this case when we rotate the board,
the position class changes, but only
among the 4 with single peg representatives in the holes shown in Figure~\ref{fig3}b.
Another way to look at Figure~\ref{fig3}b is that we can begin with $(\xstart,\ystart)$
at any ``3" or ``4", and finish at any ``3" or ``4",
if we allow peg solitaire jumps plus rotations and reflections of the board.
The third ``Type C" problems are shown in Figure~\ref{fig3}c.
The three problem types are in a sense completely separate---it
is never possible to move from a SVSS problem of one type to
another, even if you are allowed to rotate or reflect the board.

\begin{figure}[htb]
\centering
\epsfig{file=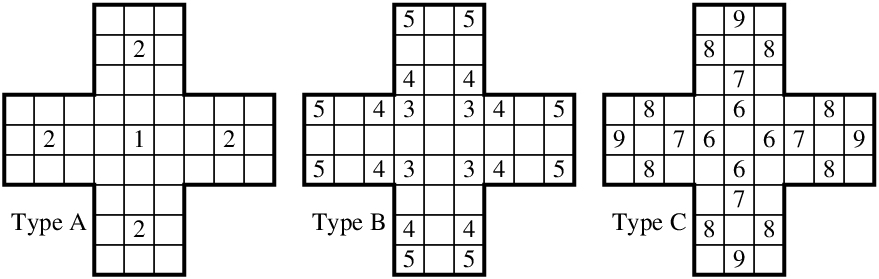}
\caption{The three types of SVSS problem on Wiegleb's 45-hole board.}
\label{fig4}
\end{figure}

\begin{figure}[htb]
\centering
\epsfig{file=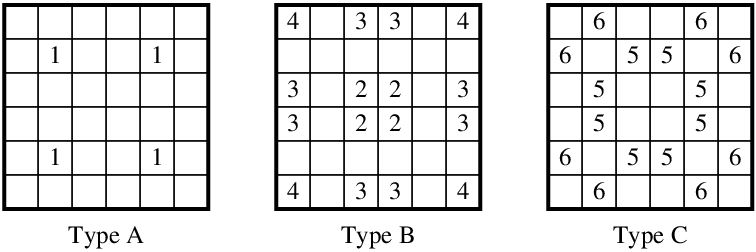}
\caption{The three types of SVSS problem on the square 6x6 board.}
\label{fig5}
\end{figure}

Wiegleb's board (Figure~\ref{fig4}) is an extension of the standard
33-hole board and has 45 holes \cite[p. 199--201]{Beasley}.
Wiegleb's board also has three problem types,
but there are more SVSS problems possible (36 in all, see \cite{BBWiegleb}).
Figure~\ref{fig5} shows the $6\times 6$ square board, which also
has three problem types.
In general, any square-lattice board will have three problem types.

\begin{figure}[htb]
\centering
\epsfig{file=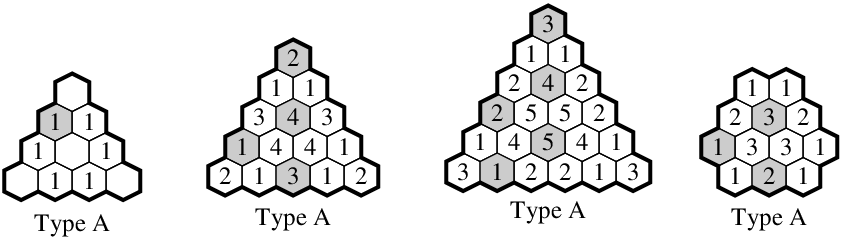}
\caption{SVSS problems on 10, 15, 21 hole triangle boards, and the 12-hole truncated triangle board.
The shaded holes show starting vacancies all in the same position class.}
\label{fig6}
\end{figure}

In triangular peg solitaire the situation is somewhat different,
for one thing there is only one type of problem, Type A.
On the triangle board of side $n$, if $n \equiv 1 \pmod 3$,
the board is not null-class and no complement problem can be solved
(see Bell \cite{BellSol}).
The only SVSS problem on the 10-hole triangle board that is solvable
is of the form: vacate $(0,1)$ finish at $(1,1)$.
This gives only one type of problem that is solvable (see Figure~\ref{fig6}a),
and any problem starting from an unmarked hole in Figure~\ref{fig6}a cannot
be solved to one peg.

The 15-hole triangle board is null-class (Figure~\ref{fig6}b).
Here there are 4 different starting locations,
the ``standard starting vacancies" for each of the 4 are shaded.
Again all problems are of type A, in a sense
all problems on this board are interrelated
(we shall see that they all must be considered together).

The 21-hole triangle board is also null-class (Figure~\ref{fig6}c),
and contains 5 different starting locations.
On this board it is possible to start with any peg missing,
and finish at any hole using peg solitaire jumps,
\textit{plus} rotating and flipping the board.
The final board in Figure~\ref{fig6}d is simply the 15-hole triangle
board with the three corner holes removed.
This ``12-hole truncated triangle" board is a good test case because all
SVSS problems are solvable on it.
It is, in fact, the smallest rotationally symmetric board with this property.

Here are several peg solitaire problems we are interested in solving:
\label{lst:problems}
\begin{packed_enumerate}
\item The \textbf{complement problem}: From a starting vacancy $(\xstart,\ystart)$,
execute an arbitrary number of jumps, then determine if the board position
can be reduced to one peg at $(\xstart,\ystart)$.
\item The \textbf{general SVSS problem}: From a starting vacancy $(\xstart,\ystart)$,
execute an arbitrary number of jumps, then determine if the board position
can be reduced to one peg (anywhere on the board).
\item The \textbf{general problem}: Given a configuration of pegs,
determine if it can be reduced to one peg (anywhere on the board).
\end{packed_enumerate}
We will consider primarily the first two problems in this paper.
We also want to solve these problems quickly---ideally
within a web browser.
In determining if a jump leads to a dead end or not,
a delay of one second is unacceptable.

It is important to realize that problems \#2 and \#3 are different.
For example, on the 33-hole board, a popular problem to solve
is ``cross" (Figure~\ref{fig1}b).
This board position is solvable to one peg,
but can never appear in the solution to any SVSS problem.
How do we know this?
Because the complement of this board position cannot be
reduced to one peg.
See Bell \cite{GPJ04} to clarify why SVSS board positions must have this property.

A subtle point is that problems \#1 and \#2 are intimately associated with the shape of the board,
while problem \#3 does not have to be associated with any board---in
the most general context we could consider problem \#3 on an infinite board.
In this sense we can see that problem \#3 is different and more general.
A fair question is, why not go for the most general and
difficult problem \#3?
The reason is that the complement problem \#1 and general SVSS problem \#2
are significantly easier, because we can take advantage of
special properties of their solutions.

\subsection{Computer solving techniques}
\label{sec:solving}

The simplest technique for solving a peg solitaire problem on a computer
is to store the sequence of completed jumps,
together with the current board position.
One then performs a depth-first search by jump
(extending the jump sequence and backtracking when
no further jump is possible).
The 15-hole triangle board can be easily solved using this
technique, but it is much slower on the 33-hole board.
The reason is that there are a large number of jump sequences
that result in the same board position, so there
is a tremendous amount of duplicated work.
This difficulty can be removed by storing board positions
seen previously in a hash table or binary tree.

An improved technique is to stop recording the jump sequence,
and look at the whole problem as a sequence of board positions (see \cite{Beasley}).
Given a set of board positions $\example$, we denote by $\descendants(\example)$ the
set of board positions that can be obtained by performing every
possible jump to every element of $\example$.
We call $\descendants(\example)$ the \textbf{descendants} of $\example$.
As a programming task calculating $\descendants(\example)$ is straightforward.
For example the set $\example$ can be stored on the disk as a
sequence of integer codes, we read each code and convert it
to a board position.
From this board position we execute every possible jump,
resulting in a large number of board positions which are stored
in a binary tree (or hash table) to remove duplicates.
This binary tree can be dumped to a file
as a sequence of codes, the set $\descendants(\example)$.

The problem that eventually occurs is that the binary tree
becomes too large to fit into memory\footnote{On Wiegleb's 45-hole board, after 22 jumps the tree contains over 2 billion elements.}.
At this point the problem is easily split into $p$ smaller pieces that
are calculated separately (or in parallel).
Let $p$ be a prime number,
chosen so that $1/p$ of the binary tree fits into memory.
We now go through the board positions in $\example$ as before,
but instead of storing each descendant in a binary tree,
we convert it to a code and write that
code to a temporary file numbered $\mbox{code}(b) \% p = \mbox{code}(b) \pmod p$.
Here the ``\%" operator represents the remainder upon
division by $p$ (as in C).
After all board positions in $\example$ have been considered,
we now go back through each of the $p$ temporary files,
filling a binary tree for each $p$ to remove duplicates and
writing the unique board positions to the disk.

In most cases we do not want to store two board positions that are
in the same symmetry class.
For example, for the central game (Figure~\ref{fig1}a),
there are four first jumps,
but these result in identical board positions which are
rotations of one another.
A convenient way to select a single representative from the symmetry
class is to use the one with the smallest code.
In the above algorithms, we use $\mbox{mincode}(b)$
in place of $\mbox{code}(b)$.

In what follows, we will denote sets of board positions by capital letters.
For a set of board positions $\example$,
we denote by $|\example|$ the number of elements in the set,
and $\overline{\example}$ is the set of complemented board positions.
In other words $b\in \overline{\example}$ if and only if $\overline{b}\in \example$.
Often these sets will contain only board positions with the same number of pegs,
so we adopt the convention that a subscript (when present) must be the number of pegs.
$\example_n$, by this convention, contains \textit{only} board positions with $n$ pegs.
We can therefore deduce that all boards in $\descendants(\example_n)$ have $n-1$ pegs and all boards in
$\overline{\example_n}$ have $N-n$ pegs
(recall that $N$ is the size of the board).
This convention is useful in understanding these sets,
for example we can immediately conclude that $\example_n\cap \example_m=\myemptyset$ when $n\ne m$.

We can also consider playing the game ``backwards", which in this
notation looks like this: let $b$ be a board position with $n$ pegs, and $B_n=\{b\}$.
Then $\overline{\descendants(\overline{B_n})}$
contains all board positions where $b$ can be reached by executing one jump.
We note that all elements of $\overline{B_n}$ have $N-n$ pegs,
$\descendants(\overline{B_n})$ have $N-n-1$ pegs, and $\overline{\descendants(\overline{B_n})}$
have $N-(N-n-1) = n+1$ pegs, as expected.

Let $\boardstart$ be the starting board position with one peg missing,
and let $\boardend$ be the target board position with one peg.
Let $\forward_{N-1}=\{\boardstart\}$ and $\backward_1=\{\boardend\}$.
We then define:
\begin{equation}
\forward_n = \descendants(\forward_{n+1}), n=N-2, N-3, \ldots, 1
\label{eq:forw}
\end{equation}
Note that, as required of our notation, every element of $\forward_n$ has exactly $n$ pegs.
This produces a ``playing forward" sequence of sets
$\forward_{N-2}, \forward_{N-3}, \ldots, \forward_1$.
We are calculating the nodes in the ``game tree", but have
lost all information about the links connecting them
(however, this link information is easily recovered).
The problem has a solution if and only if $\boardend\in \forward_1$.

A sequence of sets can also be obtained from the
finishing board position $\boardend$ by ``playing backwards" from $\backward_1=\{\boardend\}$,
\begin{equation}
\backward_n = \overline{\descendants\left(\overline{\backward_{n-1}}\right)}, n=2,3, \ldots, N-1
\label{eq:back}
\end{equation}
Again, every element of $\backward_n$ has exactly $n$ pegs,
and the problem has a solution if and only if $\boardstart\in \backward_{N-1}$.
It is worth noting that the sets $\backward_n$ contain
\textit{every} board position which can be reduced to $\boardend$ on this board.
Thus, the sets $\backward_n$ can be used to solve any problem \#3
(p. 6) which finishes at $\boardend$.
If we calculate $\backward_n$ over all possible one peg finishes,
we can solve any problem \#3 (on this board).

The set of ``winning board positions" with $n$ pegs is then
defined as
\begin{equation}
\winning_n = \forward_n \cap \backward_n
\label{eq:defwin}
\end{equation}
If we have \textit{any solution}, and play this solution
until reaching board $b$ with $n$ pegs, then it must
be that $b\in \winning_n$.
The sets $\winning_n$ tend to be much smaller than $\forward_n$ and $\backward_n$.
These winning sets $\winning_n$ are the nuggets of gold that we seek,
because they will enable us to quickly recognize when a jump leads to a dead end.

As a practical matter, to find $\winning_n$ it is not necessary to calculate every
$\forward_n$ and $\backward_n$ for each $n$ between $1$ and $N-1$ and perform their
intersection
(intersecting two sets with potentially billions of elements is not a trivial computation).
Suppose we can calculate the forward sets to $\forward_k$,
and the backward sets to $\backward_k$ for some $k$ between $N-1$ and $1$.
If the problem has a solution, then $\winning_k=\forward_k\cap\backward_k$ is not empty.
We then compute $\winning_{k-1}, \winning_{k-2}, \ldots, \winning_1=\{\boardend\}$ recursively using
\begin{equation}
\winning_n = \descendants(\winning_{n+1}) \cap \backward_n \mbox{ for } n=k-1,k-2, \ldots, 1,
\label{eq:wup}
\end{equation}
and $\winning_{k+1}, \winning_{k+2}, \ldots, \winning_{N-1}=\{\boardstart\}$ using
\begin{equation}
\winning_n = \forward_n \cap \overline{\descendants(\overline{\winning_{n-1}})} \mbox{ for } n=k+1,k+2, \ldots, N-1.
\label{eq:wdown}
\end{equation}
To calculate using equation (\ref{eq:wdown}), we take each element of $\winning_{n-1}$,
complement it, calculate all descendants and complement each.
This yields the set $\overline{\descendants(\overline{\winning_{n-1}})}$, and we now
save each element which is in common with $\forward_n$, giving us $\winning_n$.
The recursive calculations (\ref{eq:wup}) and (\ref{eq:wdown}) are much easier than
calculating all $\forward_n$ and $\backward_n$ because the sets $\winning_n$
tend to be orders of magnitude smaller.
The determination of $\winning_n$ using (\ref{eq:wup}) and (\ref{eq:wdown})
is considerably faster than the initial task of calculating the
sets $\forward_k$ and $\backward_k$.

Equations (\ref{eq:wup}) and (\ref{eq:wdown}) do not follow directly from (\ref{eq:defwin}),
so here we justify that they are correct.
For every board position $w_n\in W_n$, there must exist a 
sequence of jumps from the starting board position $b_s$ to $w_n$
(this is a direct consequence of the fact that $w_n\in F_n$),
\textit{and} a sequence of jumps from $w_n$ to the finishing board position $b_f$
(because $w_n\in B_n$).
Taken all together, these give the sequence of jumps in a solution.
Therefore, associated with every $w_n\in W_n$ there exists (at least one) sequence
of $N-1$ board positions,
\begin{equation}
\boardstart=w_{N-1}, w_{N-2}, \ldots, w_{n+1}, w_n, w_{n-1}, \ldots, w_2, w_1 = \boardend
\nonumber
\end{equation}
which show the state of the board as the solution is played.
For every 
element $w_k$ in this sequence, $w_k\in W_k$, $w_k \in \descendants(W_{k+1})$ and
$w_k \in \overline{\descendants(\overline{\winning_{k-1}})}$.
These last two statements are exactly what is needed to prove (\ref{eq:wup}) and (\ref{eq:wdown})
from (\ref{eq:defwin}).

Finally, we note that the sets $\forward_n$, $\backward_n$ and $\winning_n$
can be defined in two subtly different ways.
First, they can simply be sets of board positions.
If $\boardstart$ is the starting position for the central game (Figure~\ref{fig1}a),
then $\forward_{32}=\{\boardstart \}$, and $\forward_{31}$ has 4 elements
which are rotations of one another.
We will sometimes refer to these sets as ``raw $\forward_n$".
In most cases, however, we will consider
$\forward_n$, $\backward_n$ and $\winning_n$
as sets of symmetry classes.
Now the set $\forward_{31}$ only has a single element,
which can be taken as any representative of this symmetry class,
and generally we choose the one with the smallest $\mbox{code}()$.
These sets are called ``symmetry reduced" $\forward_n$.
If we refer to an unqualified $\forward_n$ or $\winning_n$
it can be assumed to be symmetry reduced.

\section{The complement problem}
\label{sec:comp}

For a complement problem, we have $\boardstart=\overline{\boardend}$, so that $B_1 = \overline{\forward_{N-1}}$,
and it follows from (\ref{eq:forw}) and (\ref{eq:back}) that
playing the game forward and backward amounts to the same thing,
\begin{equation}
\backward_{n} = \overline{\forward_{N-n}} ,
\end{equation}
\begin{equation}
\winning_n = \forward_n \cap \overline{\forward_{N-n}} ,
\label{eq:wcalc}
\end{equation}
and
\begin{equation}
\boxed{\winning_n = \overline{\winning_{N-n}}}
\label{eq:W}
\end{equation}
Equation (\ref{eq:W}) states that
\textit{the sets of winning board positions are complements of one another}.
This is a remarkable result, it is due to a fundamental symmetry
between pegs and holes\footnote{While the forward game jumps one peg over another peg into a hole,
we can consider the backward game as jumping one hole over another hole into a peg!
See ``Playing Backwards and Forwards" \cite[p. 817-8]{WinningWays}.}.
For a complement problem, we only need to store half the winning board positions.
In order to calculate $\winning_n$ the work is halved as well,
for we need only calculate the forward sets $\forward_n$ down to $k=\lfloor N/2 \rfloor$.
After performing the intersection (\ref{eq:wcalc}), the remaining $\winning_n$
are then calculated using Equation (\ref{eq:wup}),
which can be written using forward sets as
\begin{equation}
\winning_n = \descendants(\winning_{n+1}) \cap \overline{\forward_{N-n}} \mbox{ for } n=k-1,k-2, \ldots, 1.
\label{eq:wup1}
\end{equation}

Is storing winning board positions the most efficient technique?
During the start of a game, it does not seem so,
because all board positions that can be reached are winning.
Perhaps it is better to store ``losing board positions",
or positions from which a one peg finish at the starting
hole cannot be reached?

We could define the set of ``losing board positions"
with $n$ pegs as those elements of $\forward_n$ which are not in $\winning_n$.
It would be more efficient to store \textit{only} losing
board positions which are one jump away
from a winning board position.
Thus, we define
\begin{equation}
L_n = \descendants(\winning_{n+1}) - \winning_n
\label{Ldef}
\end{equation}
Table~\ref{table1} shows the sizes of $\forward_n$, $\winning_n$ and $L_n$ for
the 15 and 21-hole triangle boards.
All winning board positions for any corner complement problem
can be identified by storing just $95$ board positions (15-hole board)
or $26,401$ board positions (21-hole board).
If we store losing board positions as defined by Equation (\ref{Ldef}),
we need to store more than four times as many board positions.

\begin{table}[htb]
\begin{center} 
\begin{tabular}{ c | r  r  r  r | r  r  r }
 & \multicolumn{4}{|c|}{15-hole triangle board} & \multicolumn{3}{|c}{21-hole triangle board} \\
 & Raw & \multicolumn{3}{|c|}{Symmetry reduced} & \multicolumn{3}{|c}{Symmetry reduced} \\
$n$ (pegs) & \multicolumn{1}{|c|}{$|\forward_n|$} & $|\forward_n|$ & $|\winning_n|$ & $|L_n|$ & $|\forward_n|$ & $|\winning_n|$ & $|L_n|$ \\
\hline
20 & & & & & 1 & 1 & 0 \\
19 & & & & & 1 & 1 & 0 \\
18 & & & & & 4 & 4 & 0 \\
17 & & & & & 23 & 23 & 0 \\
16 & & & & & 117 & 117 & 0 \\
15 & & & & & 522 & 503 & 19 \\
14 & 1 & 1 & 1 & 0 & 1,881 & 1,690 & 185 \\
13 & 2 & 1 & 1 & 0 & 5,286 & 4,328 & 907 \\
12 & 8 & 4 & 2 & 2 & 11,754 & 8,229 & 3,288 \\
11 & 35 & 19 & 9 & 4 & 20,860 & 11,506 & 8,478 \\
10 & 122 & 62 & 18 & 20 & 28,697 & 11,506 & 14,701 \\
9 & 293 & 149 & 29 & 43 & 29,784 & 8,229 & 16,856 \\
8 & 530 & 268 & 35 & 86 & 23,263 & 4,328 & 13,063 \\
7 & 679 & 344 & 35 & 94 & 14,039 & 1,690 & 7,267 \\
6 & 623 & 317 & 29 & 89  & 6,683 & 503 & 3,005 \\
5 & 414 & 215 & 18 & 49 & 2,545 & 117 & 935 \\
4 & 212 & 112 & 9 & 29 & 774 & 23 & 211 \\
3 & 75 & 39 & 2 & 7 & 168 & 4 & 34 \\
2 & 18 & 10 & 1 & 1 & 28 & 1 & 4 \\
1 & 4 & 3 & 1 & 1 & 5 & 1 & 1 \\ \hline
Total & 3,016 & 1,544 & $\dagger$ 95 & 425 & 146,434 & $\dagger$ 26,401 & 68,954 \\
\end{tabular}
\caption{Size of $\forward_n$, $\winning_n$ and $L_n$ for the corner complement problem on the
15 and 21-hole triangle boards.
$\dagger$ Only half of the $\winning_n$ need to be stored, due to Equation~(\ref{eq:W}).} 
\label{table1}
\end{center} 
\end{table}

For the corner complement problem on the 15-hole triangle board,
the winning board positions are $\winning_n^2$.
Note that the subscript refers to the number of pegs, while the superscript
refers to the number assigned this starting vacancy in Figure~\ref{fig6}.
The sets $\winning_n^2$ give us a simple technique for determining if
we are ``on track" to solve the corner complement problem.
Suppose the current board state is $b$.
\begin{packed_enumerate}
\item If $b$ contains more than $\lfloor N/2 \rfloor=7$ pegs, complement the board position.
The board position now has $n$ pegs where $1 \le n \le \lfloor N/2 \rfloor$.
\item Calculate $\mbox{mincode}(b)$.
\item If $\mbox{mincode}(b)\in \winning_n^2$, then the complement problem can be solved
from the current board position, otherwise it cannot.
Notice that $\mbox{mincode}(b)\not\in \winning_n^2$ \textit{does not} necessarily imply that the
board position cannot be solved to one peg, just not to one peg
\textit{at the starting location}.
\end{packed_enumerate}

We also give the above algorithm in pseudocode:
\begin{footnotesize}
\begin{verbatim}
W[1] = {1}  ! The set W_1^2
W[2] = {10} ! The set W_2^2
W[3] = {28, 112} !    W_3^2
W[4] = {23, 58, 85, 120, 1108, 1616, 2076, 2210, 2272}
W[5] = {31, 93, 115, 601, 1054, 1138, 1140, 1562, 1648, 2183, 2218, 2245, 2280, 2348,
        2472, 2616, 2728, 2819}
W[6] = { 125,  633, 1086, 1111, 1594, 1621, 2191, 2253, 2275, 2289, 2343, 2467, 2589,
        2723, 2785, 2841, 2889, 3126, 3250, 3298, 3428, 3634, 3845, 4220, 4270, 4282,
        4691, 4728, 4817}
W[7] = {1567, 1651, 2235, 2365, 2413, 2537, 2731, 2793, 3159, 3196, 3320, 3374, 3388,
        3607, 3642, 3667, 3669, 3704, 3859, 3921, 4215, 4339, 4341, 4469, 4701, 4849,
        5302, 5350, 5746, 5810, 6881, 6985, 10053, 10065, 12065}

! board is the current board position
! side is the triangle board side (must be 5)
onePegFinishPossible(board, side) {
  n = CountPegs(board)
  totholes = side*(side+1)/2 ! should evaluate to 15
  if (n > totholes/2) then {
    mincode = 2^totholes - 1 - GetMaxCode(board)
    n = tot - CountPegs(board)
  }
  else mincode = GetMinCode(board)

  for (j=0; j<Size(W[n]); j++) {
    if (mincode==W[n][j]) return true
  }

  return false
}
\end{verbatim}
\end{footnotesize}

On boards with less than around 25 holes,
this test can easily be executed in a browser.
For example, when the user mouses over a peg, we can test out the jumps
from this peg and report whether the jump is ``good" or ``bad",
namely leads to a winning or losing board position.
In the web tool I have created \cite{BellNL}, the bad jumps are
humorously indicated by turning a peg into a bomb.

We have also calculated $\winning_n^1$ for the central game
on the standard 33-hole board (Figure~\ref{fig3}a).
Table~\ref{table2} shows the size of $\forward_n^1$ and $\winning_n^1$
for the central game on the 33-hole board.
These sets are large enough that the array search in the algorithm
\texttt{onePegFinishPossible()} is too slow,
and must be replaced by a faster search algorithm
for good real-time performance\footnote{In the attachments to this paper,
the sets $\winning_n^1$ are given sorted.  A good search technique is
then a binary search of a sorted array.}.
The set of $839,536$ board positions, stored in 4-byte integers,
requires $3.2$ Megabytes of memory.
Table~\ref{table3} shows results for the central game on Wiegleb's board
(Figure~\ref{fig4}a).
The set of $89,558,705$ board positions $\winning_1$ to $\winning_{22}$,
stored in two 4-byte integers, requires $680$ Megabytes of memory.

\begin{table}[htb]
\begin{center} 
\begin{tabular}{ c  r  r | c  r  r }
$n$ (pegs) & $|\forward_n|$ & $|\winning_n|$ & $n$ (pegs) & $|\forward_n|$ & $|\winning_n|$\\
\hline
32 & 1 & 1 & 16 & 3,312,423 & 230,230 \\
31 & 1 & 1 & 15 & 3,626,632 & 204,992 \\
30 & 2 & 2 & 14 & 3,413,313 & 162,319 \\
29 & 8 & 8 & 13 & 2,765,623 & 112,788 \\
28 & 39 & 38 & 12 & 1,930,324 & 68,326 \\
27 & 171 & 164 & 11 & 1,160,977 & 35,749 \\
26 & 719 & 635 & 10 & 600,372 & 16,020 \\
25 & 2,757 & 2,089 & 9 & 265,865 & 6,174 \\
24 & 9,751 & 6,174 & 8 & 100,565 & 2,089 \\
23 & 31,312 & 16,020 & 7 & 32,250 & 635 \\
22 & 89,927 & 35,749 & 6 & 8,688 & 164 \\
21 & 229,614 & 68,326 & 5 & 1,917 & 38 \\
20 & 517,854 & 112,788 & 4 & 348 & 8 \\ 
19 & 1,022,224 & 162,319 & 3 & 50 & 2 \\
18 & 1,753,737 & 204,992 & 2 & 7 & 1 \\
17 & 2,598,215 & 230,230 & 1 & 2 & 1 \\ \hline
Total & & & & 23,475,688 & $\dagger$ 839,536 \\
\end{tabular}
\caption{Size of $\forward_n$ and $\winning_n$ for the central game on the 33-hole cross-shaped board.
$\dagger$ Only half of the $\winning_n$ need to be stored, due to Equation~(\ref{eq:W}).} 
\label{table2}
\end{center} 
\end{table}

\begin{table}[htb]
\begin{center} 
\begin{tabular}{ c  r  r | c  r  r }
$n$ (pegs) & $|\forward_n|$ & $|\winning_n|$ & $n$ (pegs) & $|\forward_n|$ & $|\winning_n|$\\
\hline
43 & 1 & 1 & 32 & 3,702,227 & 348,705 \\
42 & 3 & 3 & 31 & 10,160,129 & 741,102 \\
41 & 11 & 10 & 30 & 25,647,378 & 1,483,185 \\
40 & 60 & 54 & 29 & 59,620,492 & 2,788,600 \\
39 & 297 & 236 & 28 & 127,737,457 & 4,898,948 \\
38 & 1,427 & 900 & 27 & 252,239,569 & 7,981,238 \\
37 & 6,459 & 3,007 & 26 & 458,623,402 & 11,958,747 \\
36 & 27,317 & 9,056 & 25 & 766,145,054 & 16,344,138 \\
35 & 106,347 & 24,990 & 24 & 1,172,139,707 & 20,224,817 \\
34 & 379,537 & 64,182 & 23 & 1,635,783,432 & 22,532,441 \\
33 & 1,238,520 & 154,345 & 22 & 2,073,430,928 & 22,532,441 \\ \hline
\multicolumn{3}{c|}{ } & Total & 6,586,989,754 & 89,558,705 \\
\end{tabular}
\caption{Size of $\forward_n$ and $\winning_n$ for the central game on Wiegleb's board (Figure~\ref{fig4}a).
Some elements of $\forward_n$ that cannot appear in $\winning_n$
have been removed by use of a resource count \cite{Beasley}.}
\label{table3}
\end{center} 
\end{table}

We can now solve all complement problems on the 15-hole triangle board
by calculating all $\winning_n^i$.
If we do this, we discover two problems with this technique.
The first is a degeneracy of the finishing hole
with respect to the board symmetry,
while the second storage inefficiency is that the sets $\winning_n^i$
may not be disjoint for different values of $i$.

\subsection{The symmetry degeneracy}
\label{sec:str_err}

\begin{figure}[htb]
\centering
\epsfig{file=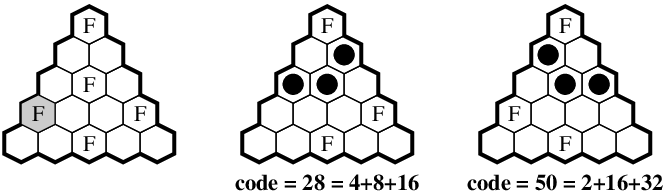}
\caption{(a) Possible finishing holes starting with $(0,3)$ empty (shaded hole).
(b) The $(0,3)$ finish can no longer be reached.
(c) This board position is in the same symmetry class as (b),
but the $(0,3)$ finish can be reached.}
\label{fig7}
\end{figure}

This problem concerns the way we have reduced the set of board positions
by using the symmetry of the board.
There is no problem in this regard to the corner vacancy,
or the central vacancy on the 33-hole standard board or Wiegleb's board.
But suppose we look at the SVSS problem on the 15-hole triangle board
starting from $(0,3)$.
According to the position class theory, the possible finishing locations
are given in Figure~\ref{fig7}a.

The problem occurs because it is possible to finish at the hole $(3,3)$,
\textit{and} that this hole is also mapped to the starting vacancy $(3,0)$
by a reflection of the board about the $y$-axis.
The board positions in Figure~\ref{fig7}b and \ref{fig7}c have the same
$\mbox{mincode}$ (28),
because they are reflections of one another,
and both can be reached starting from $(3,0)$.
The problem is that we can't finish at $(3,0)$ from Figure~\ref{fig7}b but
we can from Figure~\ref{fig7}c, yet according to our algorithm these
board positions are ``the same" (they lie in the same symmetry class).

If we create the sets $\winning_n^1$ using the symmetry reduction technique
of using $\mbox{mincode}(b)$ all will work perfectly, except
that our program will consider the finishing holes $(0,3)$ and $(3,3)$
to be the same.
One resolution of this degeneracy is to loosen our definition of
``complement problem" to include any finishing board position which
is in the same symmetry class.
In other words had we defined the problems we are trying to solve differently,
the degeneracy disappears!

But let us assume we do not want to define the problem away,
and stick with our definition of complement problem.
To resolve the degeneracy we are forced to
\textit{not to do symmetry reduction of the sets},
leaving them in their \textit{raw} state.
We then lose the simple check of calculating $\mbox{mincode}(b)$,
and checking this against the sets $\winning_n^1$.
Instead, we must figure out the symmetry transformation $S$ which takes us
from the starting board state to one peg missing at $(3,0)$.
Given any board position $b$, we then check to see if $\mbox{code}(S(b))$ is
in the set $\winning_n^1$.
Note that since we have not done symmetry reduction of the sets $\winning_n^1$,
there can be two members of this set in the same symmetry class,
so with the same $\mbox{mincode}()$.
Unfortunately, this significantly complicates our algorithm for identifying
winning board positions.

\subsection{The storage inefficiency}
\label{sec:storage}

We note from Figure~\ref{fig7} that the minimum code 28 must lie in
$\winning_3^1$, $\winning_3^2$ and $\winning_3^3$ because this mincode()
can appear during all three complement problems.
This indicates that winning board positions for different
complement problems will share members, and not just occasionally.
In fact, $\winning_n^1$ and $\winning_n^2$ have almost all of their elements in common
(we will see why soon).
This is not really a problem on the 15-hole triangle board,
because these sets are small.
It becomes more of a problem for the 21-hole triangle board,
and the 33-hole cross-shaped board.
We will discuss solutions to this problem in the next section.

\section{The general SVSS problem}
\label{sec:generalSVSS}

The key feature of complement problems which results in $\winning_n = \overline{\winning_{N-n}}$
is that the starting set $\forward_{N-1}=\{\boardstart\}$ and the finishing set
$\backward_1=\{\boardend=\overline{\boardstart}\}$ are complements of one another.
These two sets need not contain only a single board position.
For example,
let $\super{\forward}_{N-1}$ be all board positions of the same problem type with one peg missing
and $\super{\backward}_1=\overline{\super{\forward}_{N-1}}$
be all one-peg board positions of the same type.

As before we will have
$\super{\winning}_n = \super{\forward}_n \cap \super{\backward}_n = \super{\forward}_n \cap \overline{\super{\forward}_{N-n}}$ and
$\super{\winning}_n = \overline{\super{\winning}_{N-n}}$.
The winning ``superset" $\super{\winning}_n$ contains all
board positions that can appear in a SVSS problem of this type.
We see, in fact, that $\super{\winning}_n$ is the union of all $\winning_n^i$ over all
complement problems $i$ of the given type,
plus a special set which we call $\winning_n^0$
containing all board positions which can occur in SVSS
problems of this type but not in any complement problem.
We then have
\begin{equation}
\super{\winning}_n= \bigcup_{i=0,1,\ldots,p} \winning_n^i
\mbox{ where } \winning_n^0 \cap \winning_n^i = \myemptyset, i=1,2,\ldots,p
\label{eq:WW}
\end{equation}
We note that since $\super{\winning}_n = \overline{\super{\winning}_{N-n}}$ and
$\winning_n^i= \overline{\winning_{N-n}^i}$ for $i=1,2,\ldots,p$,
it must be the case that $\winning_n^0= \overline{\winning_{N-n}^0}$.

We have already seen how to calculate the complement problem
sets $\winning_n^1$, $\winning_n^2$, $\ldots$, $\winning_n^p$,
but how can we calculate $\winning_n^0$?
There really is no easy way.
We can calculate $\super{\winning}_n$ directly using all possible
starting and finishing locations,
and then subtract out each $\winning_n^i$ for $i=1,2, \ldots, p$.
This is certainly not difficult for a small board like the 15-hole triangle.
Another technique is to consider a single vacancy start, 
but any finishing hole which is of the same type.
Winning sets calculated in this manner will not satisfy (\ref{eq:W}),
but the winning boards for $n \ge \lfloor N/2 \rfloor$ will include those in $\winning_n^0$.
We then run this calculation for each starting vacancy and combine results,
finally subtracting off all complement problem boards as before.

\begin{table}[htb]
\begin{center} 
\begin{tabular}{ c | l }
$n$ (pegs) & $\winning_n^0=\super{\winning}_n$ \\
\hline
~ & \multirow{3}{*}{\includegraphics{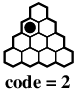}} \\
{\LARGE 1} &  \\
~ &  \\ 
~ & \multirow{3}{*}{\includegraphics{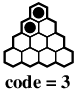}} \\
{\LARGE 2} &  \\
~ &  \\
~ & \multirow{3}{*}{\includegraphics{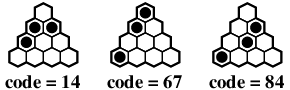}} \\
{\LARGE 3} &  \\
~ &  \\
~ & \multirow{3}{*}{\includegraphics{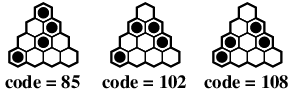}} \\
{\LARGE 4} &  \\
~ &  \\
{\large ~} & \multirow{3}{*}{\includegraphics{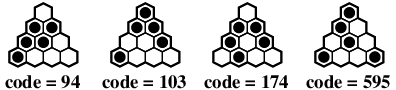}} \\
{\LARGE 5} &  \\
{\large ~} &  \\
\end{tabular}
\caption{The sets $\winning_n^0=\super{\winning}_n$ for the
10-hole triangle board.
All sets are symmetry reduced.} 
\label{table4}
\end{center} 
\end{table}

The 10-hole triangle board provides
a simple example of the sets $\winning_n^0$.
This board is not null-class,
so no complement problem is solvable.
What this means is that \textbf{all} sets $\winning_n^i$ are
empty for $i>0$ (or that the number of solvable complement problems $p=0$),
so that the only sets around are $\winning_n^0=\super{\winning}_n$.
These sets can be calculated quite easily (even by hand),
and are shown in Table~\ref{table4}.
Here we also see a board where the total number of holes $N=10$ is even,
so that $\winning_5^0=\overline{\winning_5^0}$ is equal to its own complement.
We can see that this is in fact the case, as the four board positions
in $\winning_5^0$ are listed in pairs that are complements of one other.
We must be careful to interpret the board positions in $\winning_5^0$ as symmetry classes,
the complement of the board position with code $94$ has code $2^{10}-1-94=929$,
a board position in the same symmetry class as the board with code $103$.
We conclude from Table~\ref{table4} that only $10$ essentially different
board positions (or their complement)
can appear during a solution to any SVSS problem on this board
(we only need half of the set $\winning_5^0$).

\section{The complement plus general SVSS problem}
\label{sec:compplus}

We now show how to solve either the complement problem (\#1),
or the general SVSS problem (\#2),
with very little additional storage over the general SVSS problem.
Equation (\ref{eq:WW}) shows that these two problems are closely related.
The only difficulty involves storing board positions efficiently (without duplicates),
and dealing with the degeneracies introduced in section~\ref{sec:str_err}.
In this section we use as an example the 15-hole triangle board.

We deal first with the storage problem.
How can we store all board positions in $\winning_n^i$ without duplication?
Here we are considering all the problems of a certain type on a board,
and $i=1,2,\ldots,p$ ranges over the total number of problems of this type.
On the 15-hole triangle board, there is only one type
(shown in Figure~\ref{fig6}b) with $p=4$ different problems.
The obvious solution is to take all possible combinations
of the 4 problems, $2^p=16$ possibilities,
and for each combination we store all boards
common to this combination of problems.
We can think of these sets $W_n$ with a superscript given not by the
problem number $i$, but the \textit{index} of the combination
of complement problems that this board position can occur in
(ranging from $0$ to $16$).
Index $0$ remains the same as problem $0$:
$\winning_n^0=\winning_n^{index=0}$ is the set of $n$-peg
board positions that can occur during a SVSS problem,
but not in any complement problem.

\begin{table}[phtb]
\begin{center} 
\begin{tabular}{ c | p{0.8cm}  p{0.8cm}  p{0.8cm}  p{0.8cm}  p{0.8cm}  p{0.8cm}  p{0.8cm} | p{1.0cm} }
 & \multicolumn{7}{|c|}{$|\winning_n^{index}|$ for $index=0$ to $7$} & \\
Pegs ($n$) & \rtj{0} & \rtj{1} & \rtj{2} & \rtj{3} & \rtj{4} & \rtj{6} & \rtj{7} & Total \\
\hline
1 & \rtj{1} & \rtj{1} & \rtj{1} & \rtj{0} & \rtj{1} & \rtj{0} & \rtj{0} & \rtjl{4} \\
2 & \rtj{1} & \rtj{0} & \rtj{0} & \rtj{1} & \rtj{2} & \rtj{0} & \rtj{0} & \rtjl{4} \\
3 & \rtj{3} & \rtj{0} & \rtj{0} & \rtj{1} & \rtj{7} & \rtj{0} & \rtj{1} & \rtjl{12} \\
4 & \rtj{5} & \rtj{0} & \rtj{0} & \rtj{4} & \rtj{19} & \rtj{1} & \rtj{4} & \rtjl{33} \\
5 & \rtj{10} & \rtj{0} & \rtj{1} & \rtj{8} & \rtj{49} & \rtj{4} & \rtj{8} & \rtjl{80} \\
6 & \rtj{7} & \rtj{0} & \rtj{2} & \rtj{11} & \rtj{93} & \rtj{6} & \rtj{13} & \rtjl{132} \\
7 & \rtj{4} & \rtj{0} & \rtj{2} & \rtj{12} & \rtj{129} & \rtj{7} & \rtj{18} & \rtjl{172} \\ \hline
Total & \rtj{31} & \rtj{1} & \rtj{6} & \rtj{37} & \rtj{300} & \rtj{18} & \rtj{44} & \rtjl{437} \\
\end{tabular}
\caption{A count of winning board positions $\winning_n^{index}$ on the 15-hole triangle board.
\textit{index} ranges from $0$ to $15$ but $5$ and those beyond $7$ are all empty.} 
\label{table5}
\end{center} 
\end{table}

For example, since $index=7=0111$ in binary, then
$\winning_n^{index=7}$ contains all $n$-peg positions that are common to
problems 1, 2 and 3.
We note the since problem 4 is unsolvable as a complement problem \cite{BellSol},
$\winning_n^4=\myemptyset$ and all sets $\winning_n^{index}$ with $index$
between $8$ and $15$ are also empty.

A trickier question is how to resolve the degeneracy at the $(0,3)$
starting location.
The sets $\winning_n^1$ cannot be symmetry reduced,
yet $\winning_n^2$ and $\winning_n^3$ are symmetry reduced,
and $\winning_n^4=\myemptyset$.
We see from Figure~\ref{fig7} that $\winning_n^1$ contains code $50$,
while $\winning_n^2$ contains mincode $28$,
a board position in the same symmetry class.
The solution is to use an algorithm which keeps all codes in $\winning_n^1$
but removes all symmetry equivalents in the intersecting sets.
This is the reason that the degenerate starting locations
are indexed first.

Table~\ref{table5} summarizes the number of board positions by $index$ and number of pegs $n$.
The total number of board positions over all sets is $437$, which is only $10$ more than
are needed to solve the general SVSS problem by itself.
There are $10$ board positions in $\winning_n^1$ which are reflections of other board positions
in $\winning_n^1$ (an example which is useful to check by hand are the board positions with
codes 93 and 563, these 5-peg board positions end up in \textit{index} 7).
The first row of Table~\ref{table5} tells us which complement problems are solvable,
namely those with \textit{index} 1 (problem 1, the $(0,3)$ complement),
2 (problem 2, the $(0,0)$ complement) and 4 (problem 3, the $(2,4)$ complement).
Again, the fact that $W_1^8=\myemptyset$ is equivalent to the statement that problem 4,
the $(1,2)$ complement, is not solvable.

\begin{figure}[htb]
\centering
\epsfig{file=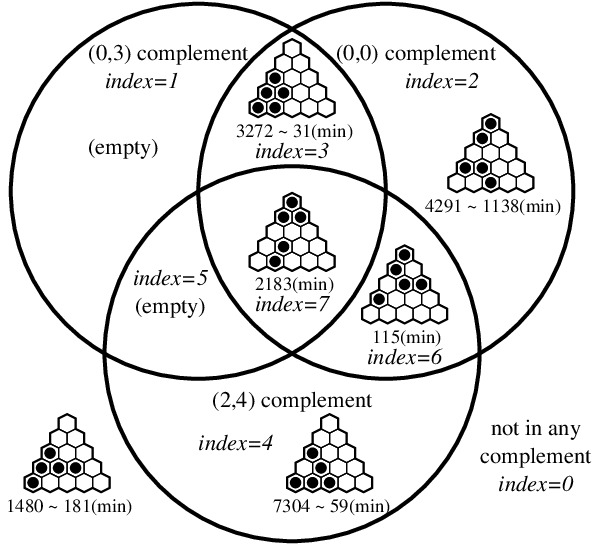}
\caption{A Venn diagram showing sample boards with 5 pegs by index.}
\label{fig8}
\end{figure}

Figure~\ref{fig8} shows representative 5-peg board positions
in $\winning_5^{index}$ for values of the \textit{index} 0 to 7.
Let us interpret two of the board positions in this diagram.
The board position with code $1480$ is in \textit{index} 0, meaning
that this board position cannot appear in any complement problem.
This board position has a mincode of $181$.
We can finish with one peg from this board position at $(0,0)$ or
$(0,3)$, but we cannot start from either of these holes and reach
this board position.
But it must be possible to reach this board position from some start,
and it turns out this start is $(2,4)$.

The board position with \textit{index} 7 is $2183$, which is the mincode.
We can play from this board position to finish at $(0,0)$, $(0,3)$,
$(3,3)$ or $(2,4)$, and we can also reach this board position starting
from $(0,0)$, $(0,3)$ or $(2,4)$.
Therefore, this board position can be reached during any of the
three (solvable) complement problems.

The board position with \textit{index} 6 is $115$, which is the mincode.
We can play from this board position to finish at $(0,0)$, $(0,3)$
or $(2,4)$, and we can reach this board position from $(0,0)$ or $(2,4)$.
Therefore, this board position can be reached during the solution to
the $(0,0)$ or $(2,4)$ complements, so is in \textit{index} 6.

We can see in Figure~\ref{fig8} (and Table~\ref{table5}) that the \textit{index} 1
and \textit{index} 5 sets are empty,
why is that?
The reason is that $W_1^2$ contains a single board position with mincode 10 and
pegs at $(0,1)$ and $(0,2)$.
This means that all solutions to the $(0,3)$ complement must \textit{begin and finish}
with the jump from $(0,1)$.
Therefore, any board position which can be reached during the $(0,3)$ complement
(except for the starting board position) can be reached during the $(0,0)$
complement\footnote{The converse is not true, do you see why?
The $(0,0)$ complement can begin with a jump from $(0,2)$,
but end with a jump from $(2,2)$.}.
This is exactly what it means for the \textit{index} 1 and \textit{index} 5 sets to be empty,
except for the one peg starting board position for the $(0,3)$ complement.

We now present pseudocode for identification of winning board positions for either
the complement problem (\#1) or the general SVSS problem (\#2):

\begin{footnotesize}
\begin{verbatim}
W[5][1] = {16,64,1,8} ! W^1, index=0,1, ... 15, for the 15-hole triangle board
End[5][1] = {1,2,3,3,4,4,4,4,4,4,4,4,4,4,4,4} ! Ends of each index 0-15
... see TriangleWinning/Triangle5ByIndex.txt ...
! board is the current board position
! side is the triangle board side (4,5, or 6)
! i is the number of this SVSS problem
! ksym is the symmetry code of the starting board position
! comp is true for complement problems, otherwise any finish is assumed
problemIsSolvable(board, side, i, ksym, comp) {
  int code[6]

  totholes = side*(side+1)/2

  if (side==4) { ! 10-hole triangle board
    topindex = 1
    degen = 0
  }
  if (side==5) { ! 15-hole triangle board
    topindex = 2^4 ! power set of 4 problems
    degen = 1 ! number of degenerate problems
  }
  if (side==6) { ! 21-hole triangle board
    topindex = 2^5 ! power set of 5 problems
    degen = 2 ! number of degenerate problems
  }

  n = CountPegs(board)
  if (n > tot/2) then {
    code[0] = 2^tot - 1 - Code(board)
    n = tot - CountPegs(board)
  }
  else code[0] = Code(board)

  ! Get the 6 symmetry codes, code[0] to code[5]
  code[1] = rotatecode(code[0])
  code[2] = rotatecode(code[1])
  code[3] = reflectcode(code[2])
  code[4] = rotatecode(code[3])
  code[5] = rotatecode(code[4])
 
  if (comp) { ! complement problem
    kStart = 0
    kEnd = 6
    if (i<=degen) {
      kStart = kSym
      kEnd = kSym + 1
    }
    for (index=1; index<topindex; index++) {
      if ((1<<i) & index) { ! true if the i'th bit of index is set
        for (j=End[side][n][index-1]; j<End[side][n][index]; j++) {
          for (k=kStart; k<kEnd; k++) if (code[k]==W[side][n][j]) return true
        }
      }
    }
  }
  else { ! comp=false, finish anywhere
    for (j=0; j<End[side][n][topindex]; j++) {
       for (k=0; k<6; k++) if (code[k]==W[side][n][j]) return true
    }
  }

  return false
}
\end{verbatim}
\end{footnotesize}

\section{Results from winning board calculations}
\label{sec:results}

\subsection{How badly can you play?}
\label{sec:badly}

After calculating the winning sets $W_n$ we can determine
the first possible ``dead end" for a complement problem.
This is the shortest jump sequence after which the goal can no longer be reached.
A well-known sequence of 4 jumps is a dead end for the central game (proved in \cite[p. 115]{Beasley}).
We can find this board position computationally by looking for the first set where
$F_n \ne W_n$, from Table~\ref{table2} we can see that this happens after 4 jumps at $n=28$ pegs.

\begin{figure}[htb]
\centering
\epsfig{file=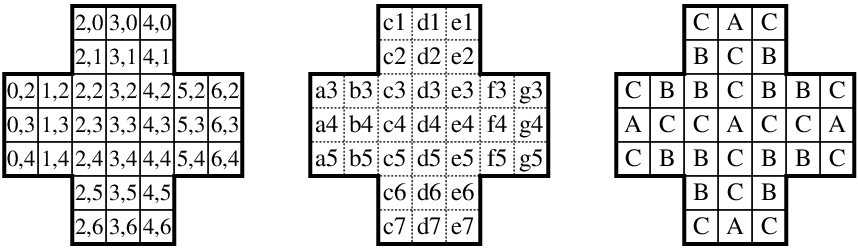}
\caption{For the cross-shaped 33-hole board: hole coordinates, Beasley's notation,
and the three types of problem (from Figure~\ref{fig3}).}
\label{fig_notation}
\end{figure}

We have calculated the first dead end for all seven complement problems, the
results are shown in Table~\ref{table6}.
Surprisingly, there are three complement problems that can be lost in only three jumps.
All complement problems also have a unique jump sequence leading to a dead end
(up to symmetry and jump order), with the exception of the c3-complement.
Jumps are presented in the alphanumeric format used by John Beasley \cite{Beasley}
(Figure~\ref{fig_notation}).

\begin{table}[htb]
\begin{center} 
\begin{tabular}{ c l | r | c | l | l }
& complement & \multirow{2}{*}{$\sum |\winning_n|$} & & & can still\\
type & problem & & jumps & first dead end & finish at\\
\hline
\multirow{2}{*}{A} & $(3,3)$ or d4 & 839,536 & 4 & {\small d2-d4, d5-d3, b4-d4, d3-d5} & nowhere \\
& $(3,0)$ or d1 & 99,982 & 3 & {\small d3-d1, b3-d3, e3-c3} & a4, d4, g4 \\ \hline
\multirow{2}{*}{B} & $(2,2)$ or c3 & 20,836,420 & 6 & {\small c1-c3, e2-c2, c3-c1, d4-d2, b4-d4, f3-d3} & nowhere \\
& $(2,1)$ or c2 & 12,372,794 & 4 & {\small c4-c2, a4-c4, d4-b4, f4-d4} & c5, f5\\ \hline
\multirow{3}{*}{C} & $(3,2)$ or d3 & 6,420,923 & 3 & {\small d1-d3, d4-d2, d6-d4} & a3, g3 \\
& $(3,1)$ or d2 & 760,164 & 3 & {\small d4-d2, b4-d4, e4-c4} & a5, d5, g5 \\
& $(2,0)$ or c1 & 13,918,925 & 5 & {\small c3-c1, e2-c2, d4-d2, c2-e2, b4-d4} & c7 \\
\end{tabular}
\caption{Shortest dead ends for all complement problems on the 33-hole board.
All sets $\winning_n$ are symmetry reduced.} 
\label{table6}
\end{center} 
\end{table}


Note that the jumps shown in Table~\ref{table6} are dead ends only for the complement problem.
As indicated in the rightmost column, in many cases it is still possible to finish with one peg,
just not at the location of the initial vacancy.
Suppose we want to find the first time we can reach
a board position from which a single peg finish is impossible?
This question can be answered by calculating the winning board sets $\super{\winning}_n$
as we did in section~\ref{sec:generalSVSS},
starting from all single vacancies of the same type.

\begin{figure}[htb]
\centering
\epsfig{file=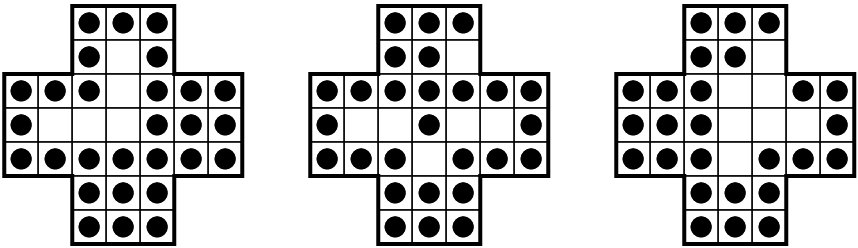}
\caption{Fastest dead board positions for problems of Type~A (starting from d1, a4 or d4, after 4 jumps),
Type~B (from c2, after 5 jumps), Type~C (from c1, c4 or f4, after 6 jumps).}
\label{fig9}
\end{figure}

To distinguish from a complement problem ``dead end",
we will now call any board position from which a single peg finish is impossible a ``dead board position".
For any Type~A single vacancy (d1, a4, d4, g4 or d7) a sequence of 4 jumps can land you in a dead board position.
For d4, these 4 jumps are given in Table~\ref{table6},
and the resulting board position is shown in Figure~\ref{fig9}a.
From the Type B vacancy at c3, we can play the same six jumps as in Table~\ref{table6} and reach a dead board position,
and many other 6 jumps sequences will also work (but none shorter).
From the Type~B vacancy at c2, there is a \textit{unique} sequence of 5 jumps which ends at a 
dead board position, shown in Figure~\ref{fig9}b.
From any Type~C vacancy, any 5 jumps can be made, and it is always possible to finish with one peg.
There are many combinations of 6 jumps which can lead to dead board positions.
Figure~\ref{fig9}c shows one board position which can be reached in 6 jumps from c1, c4 or f4.
This dead board position, or a rotation or reflection of it,
can therefore be reached from any Type C single vacancy.

At this point we have only a ``computer proof" that the board positions in
Figure~\ref{fig9}b and \ref{fig9}c cannot be reduced to a single peg.
By applying some of the techniques in Beasley's book \cite{Beasley}, it is possible to prove
this analytically (try resource counts or Conway's balance sheet \cite[p. 101-116]{Beasley}).
Of course, this does not prove that these are the \textit{first possible} dead board positions.

One can also calculate fastest dead ends and fastest dead board positions
for problems on other board types.
The central game on Wiegleb's Board can be lost in only three jumps,
this can be deduced from Table~\ref{table3} (see \cite{BellWeb} for the three jumps).
On the 15-hole triangle board, it is possible to reach a dead board position after only one
jump\footnote{Start with $(0,1)$ vacant, and jump the peg from $(2,3)$.}!

\subsection{How many wins are there?}
\label{sec:wins}

\begin{figure}[htb]
\centering
\epsfig{file=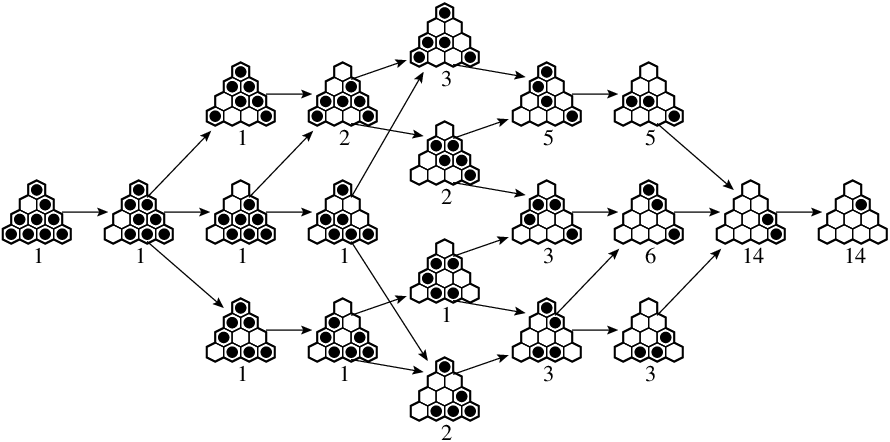}
\caption{A winning game graph for the 10-hole triangle board.}
\label{fig10}
\end{figure}

We can also use $W_n$ to count the number of solutions to any SVSS problem.
Suppose we take the set of winning board positions,
and add directed edges between board positions
related by peg solitaire jumps.
This results in a directed graph of the peg solitaire win, Figure~\ref{fig10} shows this
graph for the 10-hole triangle board.

The total number of wins is simply the total number of ways to traverse this graph.
Here is a simple algorithm for calculating this:
we will label each vertex (board position) with a number,
which will be the number of paths from the start to that board position.
We first label the starting board position with a 1.
We then consider all descendants of labeled board positions,
and label each as the sum of all incoming edges.
We continue this process until we reach the final board position.
As can be seen in Figure~\ref{fig10}, the total number of ways to traverse the graph,
or the number of solutions, is 14.

We have run this counting algorithm on the central game on the 33-hole board,
using the \href{http://arxiv.org/src/0903.3696v4/anc}{ancillary program} ``count.cpp".
In the calculation we have eliminated symmetry considerations to be sure that we calculate every solution.
This graph therefore has more than the $839,536$ nodes of the symmetry reduced winning sets.
The program calculates that the total number of winning games is $40,861,647,040,079,968 \approx 4.1\times 10^{16}$,
in perfect agreement with the figure calculated by Bill Butler \cite{DurangoBill}.
We have counted solutions for all seven complement problems, see Table~\ref{table7}.

\begin{table}[htb]
\begin{center} 
\begin{tabular}{ c c | r | r | c }
 & complement & \multicolumn{2}{|c|}{number of solutions} & calculation\\
type & problem & exact & approx. & time $\dagger$\\
\hline
\multirow{2}{*}{A} & $(3,3)$ or d4 & 40,861,647,040,079,968 & $4.1\times 10^{16}$ & 7 minutes \\
& $(3,0)$ or d1 & 841,594,661,434,808 & $8.4\times 10^{14}$ & 8 seconds \\ \hline
\multirow{2}{*}{B} & $(2,2)$ or c3 & 138,409,681,956,904,365,268 & $1.4\times 10^{20}$ & 90 minutes \\
& $(2,1)$ or c2 & 17,385,498,352,036,301,092 & $1.7\times 10^{19}$ & 22 minutes \\ \hline
\multirow{3}{*}{C} & $(3,2)$ or d3 & 8,940,989,276,947,390,168 & $8.9\times 10^{18}$ & 20 minutes \\
& $(3,1)$ or d2 & 30,997,283,487,697,056 & $3.1\times 10^{16}$ & 1 minute \\
& $(2,0)$ or c1 & 2,343,652,440,537,181,612 & $2.3\times 10^{18}$ & 25 minutes \\
\end{tabular}
\caption{The number of solutions to complement problems on the 33-hole cross-shaped board.
$\dagger$ Run time is on a $2.66$ GHz Windoze machine (single processor).}
\label{table7}
\end{center} 
\end{table}

\section{Online games and software}
\label{sec:anc}

I have created a Javascript game \cite{BellNL} for playing peg solitaire on the
10, 15 and 21-hole triangle boards,
as well as the 12 and 18-hole truncated triangle boards.
The game can begin from any starting vacancy,
and the program will point out all bad jumps which leave the set of winning board positions
(when a user hovers over a peg that can make a jump,
the GUI indicates this by turning that peg into a bomb).
The game can be specified as either a complement problem or the general problem
with a one peg finish anywhere on the board.
The algorithm used to identify winning board positions
is the algorithm \texttt{problemIsSolvable()} in this paper converted to Javascript.
This program can also find a solution from any (solvable) intermediate board position
by testing jumps chosen at random.

The \href{http://arxiv.org/src/0903.3696v4/anc}{ancillary files}
for this paper include the following directories:

\begin{itemize}
\item \texttt{Triangular1.2} -- a collection of html and javascript programs which can solve
SVSS problems on triangle boards of arbitrary size.

\item \texttt{NeverLose1.5} -- a collection of html and javascript programs which can identify
all winning positions on the 10, 12, 15, and 21-hole triangle boards
(source files for \cite{BellNL}).

\item \texttt{TriangleWinning} -- contains text files of the board positions $\super{\winning}_n$,
for triangle and truncated triangle boards.
Each subdirectory contains data for each of the 5 boards given in Table~\ref{table8}.
The data $\super{\winning}_n$ for each board comes in two versions: first 
an ``AnyFinish" version which is simply the codes $\super{\winning}_n$ all symmetry reduced
and sorted for each $n$.
If this data is used in the simple algorithm \texttt{onePegFinishPossible()}
it will point out all board positions for problem \#2 which can finish with one peg.
The second ``ByIndex" version of $\super{\winning}_n$ has the board positions sorted by index,
and degenerate starting locations are not symmetry reduced.
When applied in the algorithm \texttt{problemIsSolvable()},
we can identify when a board position can appear in any complement problem \#1 as well
as any one peg finish (problem \#2).
It is the second set of data that you will find in the Javascript web program \texttt{NeverLose1.5}.
Table~\ref{table8} contains a summary of how many board positions are included in these two
data sets, for each of the five boards.

\begin{table}[htb]
\begin{center} 
\begin{tabular}{ l | r  r }
 & \multicolumn{2}{| c}{number of codes needed to solve} \\
board & any finish (\#2) & any complement problem (\#1) \\
\hline
Triangle10 & 12~$\dagger$ & 0 \\
TruncTriangle12 & 136 & 147 \\
Triangle15 & 427 & 437 \\
TruncTriangle18 & 8,621 & 11,444 \\
Triangle21 & 76,981 & 110,647 \\
\end{tabular}
\caption{The number of codes needed to solve all problems on a particular board.
$\dagger$ These board positions are shown in Table~\ref{table4}.}
\label{table8}
\end{center} 
\end{table}

\item \texttt{pegs} -- a collection of C++ routines for calculating winning sets $\winning_n$
for the 33-hole cross-shaped board using Equation (\ref{eq:wup1}).
Table~\ref{table2} gives the sizes these sets for the central game.
If printed out in a text file, all 839,536 winning board positions for the central
game take up about 9MB.  These can be printed out by the above program,
the beginning and end of these sets is given below.

\begin{small}
\noindent
$\winning_1=\{65536=2^{16}\}$,
$\winning_2=\{528=2^4+2^9\}$,\newline
$\winning_3=\{400=2^4+2^7+2^8, 212992=2^{14}+2^{16}+2^{17}\}$,\newline
$\winning_4=\{$153, 1680, 16688, 17928, 66432, 82976, 147984, 352256$\}$,\newline
$\winning_5=\{$158, 692, 793, $\ldots$, 4554760, 6684688, 8601616$\}$,\newline
...\newline
$\winning_{14}=\{$53247, 56831, 57279, $\ldots$, 2651879594, 2655805539, 3098292302$\}$,\newline
$\winning_{15}=\{$127999, 128895, 129791, $\ldots$, 3793449102, 3793531059, 3793629859$\}$,\newline
$\winning_{16}=\{$126975, 130559, 229359, $\ldots$, 3864553651, 3928764638, 3929805043$\}$.\newline
\end{small}

The pegs directory also contains the program ``count.cpp", which counts the number of winning games
after ``pegs.cpp" has been run.

\item \texttt{old\_pegs} -- a collection of older C++ routines for calculating winning sets $\winning_n$
and $\super{\winning}_n$ for arbitrary board shapes.
They can handle triangular boards and even diagonal jumps, and calculate by move rather than by jump.
One version can handle boards with up to 48 holes, a second version can handle
boards of unlimited size.
These programs calculated the winning sets for the triangle boards.
See \texttt{readme.txt} in this directory for a complete description of these files.

\item \texttt{FigGen} -- contains the C++ program which generated the figures in this paper.
This program generates text files (.fig extension) which are input files
for the free UNIX drawing program \texttt{Xfig}.
Simply open the .fig file in \texttt{Xfig} and export to .eps to be used by \LaTeX .
\end{itemize}

\section{Summary}
\label{sec:summary}

We have introduced some techniques in peg solitaire for calculating winning board positions.
These winning board positions are useful because with them we can create a program which will
begin from any single vacancy start and identify all board positions from which
\begin{packed_enumerate}
\item it is possible to finish where we started
(the complement problem is still solvable) or
\item it is possible to finish with one peg.
\end{packed_enumerate}

We note that these techniques, by design, work \textit{only} for board positions which
can appear during SVSS problems.
Suppose we consider an \textit{arbitrary} board position $b$.
A reader may conclude from Equation (\ref{eq:W}) that if $b$ is solvable to one peg,
then $\overline{b}$ must also be solvable to one peg.
\textbf{This is false!}
Consider, for example, the board position ``cross" Figure~\ref{fig1}b.
The statement of Equation (\ref{eq:W}) applies only to board positions which can arise
during SVSS problems, not any arbitrary board position.
This boils down to the difference between solving problems \#1 and \#2 compared with problem \#3
as discussed on page~\pageref{lst:problems} (section~\ref{sec:SVSS}).

We have also shown how the winning board positions can be used to
answer some difficult questions about peg solitaire.
For the 33-hole board, we have calculated the first possible complement problem 
dead end as well as the first possible dead board position.
We have also counted the number of solutions to all complement problems.

\newpage



\begin{thebibliography}{mybib} 

\bibitem{Beasley} J. Beasley, \textit{The Ins and Outs of Peg Solitaire}, Oxford Univ. Press, 1992.

\bibitem{WinningWays} E. Berlekamp, J. Conway and R. Guy, Purging pegs properly, in \textit{Winning Ways for Your Mathematical Plays}, 2nd ed., Vol. 4, Chap. 23: 803--841, A K Peters, 2004.

\bibitem{BBWiegleb} G. Bell and J. Beasley, New problems on old solitaire boards, Board Game Studies, \textbf{8} (2006),
\href{http://arxiv.org/abs/math/0611091}{\tt http://arxiv.org/abs/math/0611091}

\bibitem{GPJ04} G. Bell, Triangular peg solitaire unlimited, {\it Games and Puzzles J.} \#36 (2004),\newline
\href{http://www.gpj.connectfree.co.uk/gpjr.htm}{\tt http://www.gpj.connectfree.co.uk/gpjr.htm}\newline
\href{http://arxiv.org/abs/0711.0486}{\tt http://arxiv.org/abs/0711.0486}

\bibitem{BellSol} G. Bell, Solving triangular peg solitaire, \textit{J. Integer Sequences}, 08.4.8, \textbf{11} (2008),\newline
\href{http://arxiv.org/abs/math/0703865}{\tt http://arxiv.org/abs/math/0703865}

\bibitem{BellNL} G. Bell, ``Never Lose" peg solitaire game,\newline
{\footnotesize\href{http://home.comcast.net/~gibell/pegsolitaire/Tools/Neverlose/Triang.htm}
{\tt http://home.comcast.net/~gibell/pegsolitaire/Tools/Neverlose/Triang.htm}}

\bibitem{BellFr} G. Bell, A fresh look at peg solitaire, \textit{Math. Mag.} \textbf{80} (2007), 16--28.

\bibitem{BellWeb} G. Bell, Peg Solitaire web site,
{\footnotesize\href{http://home.comcast.net/~gibell/pegsolitaire/}
{\tt http://home.comcast.net/~gibell/pegsolitaire/}}

\bibitem{DurangoBill} B. Butler, \href{http://www.durangobill.com/Peg33.html}{\tt http://www.durangobill.com/Peg33.html}
 
\end{thebibliography}
\end{document}